\magnification=\magstep 1
\hbadness=10000
\tolerance=10000
\hsize=6.5truein
\vsize=8.5truein
\parskip=0pt
\baselineskip=18 pt
\def\noi{\noindent}
\headline={\hss{}\hss}
\footline={\ifnum\pageno<2 \hss{}\hss\else\tenrm\hss\folio\hss\fi}

\def\section#1.#2.{\baselineskip=12pt\goodbreak\vskip1truecm\hangindent
            \hangafter=1\noindent{\bf #1.}\ \ {\bf #2}\nobreak\par\vskip 
5pt \baselineskip=20pt\ignorespaces}
\def\Def #1 {\medbreak\noi{\bf Definition #1}\ \ }
\def\Prop#1 {\medbreak\noi{\bf Proposition #1}\ \ }
\def\Lemma#1 {\medbreak\noi{\bf Lemma #1}\ \ }
\def\Cor#1 {\medbreak\noi{\bf Corollary #1}\ \ }
\def\Remark{\medbreak\noi{\bf Remark:}\ \ }

\def\vs{\vskip 12pt}
\def\proof{\noi{\bf Proof}:\ \ }
\def\squar{\vbox{\hrule\hbox{\vrule height 6pt \hskip
6pt\vrule}\hrule}}
\def\sqr{{\unskip\nobreak\hfil\penalty50\hskip2em\hbox{}\nobreak\hfil
{\squar}\parfillskip=0pt\finalhyphendemerits=0\par}}
\def\endproof{\sqr\par\vs}

\def\a{\alpha}
\def\la{\lambda}

\def\var{\varphi}
\def\br{$\bf R_+$}
\def\bz{$\bf Z_+$}
\def\gid{g.i.d.}
\def\qgid{q.g.i.d.}
\def\lqgid{$\la$-q.g.i.d.}
\def\nid{$\cal N$-i.d.}
\def\lqnid{$\la$-quasi-$\cal N$-i.d.}
\def\AR1{$AR(1)$}
\def\calgl{{\cal G}_\la}
\def\calglo{{\cal G}_{\la_1}}
\def\calglt{{\cal G}_{\la_2}}
\def\calgstar{{\cal G}_*}
\def\calginf{{\cal G}_\infty}

\newcount\refno
\refno=0
\def\ref{\advance\refno by 1\medbreak\item{\the\refno.}}

\def\bibs{\vskip 1truecm\baselineskip12pt\parskip 3pt
                         \centerline{\bf References}\par\nobreak}

\centerline{\bf Quasi-geometric infinite divisibility}

\vskip 1 cm

{\baselineskip=14 pt \centerline{Nadjib Bouzar}
\centerline{Department of Mathematics and Computer Science}
\centerline{University of Indianapolis}
\centerline{Indianapolis, IN 46227, USA}
\centerline{\underbar{e-mail}: nbouzar@uindy.edu}}

\vskip 1.5cm

\centerline{\bf Abstract}

\vskip .5cm

The object of this paper is to introduce and study the concept of quasi-geometric infinite divisibility for distributions on $\bf R_+$. These distributions arise as mixing distributions of (discrete) geometric infinitely divisible Poisson mixtures. Several characterizations and closure properties are presented. A connection between quasi-geometric infinite divisibility and log-convex (log-concave) distributions is established. A generalized notion of quasi-infinite divisibility is also discussed.

{\baselineskip=12pt
\footnote{}{\hskip -.3 in {\it Keywords and Phrases:} Poisson mixture, Laplace-Stieltjes transform, probability generating function, log-concavity, log-convexity.}

\footnote{}{\hskip -.3 in {\it AMS Subject Classification (2000)}: 60E07}}

\filbreak

\goodbreak\noi {\bf 1.  Introduction}

\smallskip

A real-valued random variable (rv) $X$ is said to have a geometrically infinitely divisible ({\gid}) distribution  if for any $p\in (0,1)$, there exits a sequence of iid, real-valued rv's $\{X_i^{(p)}\}$ such that
$$
X\buildrel d \over =\sum_{i=1}^{N_p}X_i^{(p)},\eqno(1.1)
$$
where $N_p$ has the geometric distribution
$$
P(N_p=k)=p(1-p)^{k-1}, \quad k=1,2,\cdots,\eqno(1.2)
$$
and $N_p$ and $\{X_i^{(p)}\}$ are independent. This definition is due to Klebanov et al. (1984). The authors also introduced the related concept of geometric stability. The theory of {\gid} distributions and geometrically stable distributions parallels nicely the theory of classical infinite divisibility as shown by Klebanov et al. (1984) and a number of other authors in a series of subsequent articles individually referenced in Kozubowski and Rachev (1999) (see also the monograph by Gnedenko and Korolev (1996) for generalizations). Aly and Bouzar (2000) studied the case of {\gid} distribtutions on {\bz}:$=\{0,1,2,\cdots\}$ and {\br}:$=[0,\infty)$. Random summation schemes such as (1.1) also turned out to be very useful in economics and finance (see Kozubowski and Rachev (1994, 1999) and references therein and Gnedenko and Korolev (1996)) and in queueing theory (Jacobs (1986)).

A Poisson mixture is a distribution on {\bz} that results from the mixing of a Poisson distribution by a distribution on {\br}. Poisson mixtures play an important role in distribution theory, particularly in the areas of infinite divisibility, self-decomposability, and stability of probability distributions on {\br} (see for example the monograph by Steutel (1970), Puri and Goldie (1979), van Harn and Steutel (1993), and Aly and Bouzar (2000)). Moreover, van Harn and Steutel (1993) and Pakes (1995) used Poisson mixtures to solve stability equations for {\br}-valued processes with stationary independent increments. Kebir (1997) derived several characterization theorems in renewal theory via Poisson mixtures. Poisson mixtures have also been identified as very reasonable models for a variety of random phenomena (cf. Johnson et al. (1992)).

The purpose of this paper to introduce and study the concept of quasi-geometric infinite divisibility ({\qgid}) for distributions on {\br}. Essentially, these are distributions that arise as mixing distributions of {\gid} Poisson mixtures (cf. definitions below). Our approach follows that of Puri and Goldie (1979). We present several characterizations of {\qgid} distributions. Closure properties are also obtained and various examples and counterexamples are given. We obtain necessary and sufficent conditions for a distribution on {\bz} to be a Poisson mixture generated by a {\gid} distribution. We establish a connection between the important concept of log-convexity (log-concavity) and the {\qgid} property by way of L\'evy measures. In the process, we derive a number of new characterizations of {\gid} distributions on {\br}. Finally, a generalized notion of quasi-infinite divisibility is also introduced.   
    
In the remainder of this section we recall a few useful facts that will be used throughout the paper. A distribution with support in {\bz} is {\gid} if and only if its probability generating function (pgf) $P(z)$ has the form
$$
P(z)=\bigl(1+c(1-Q(z)\bigr)^{-1},\qquad |z|\le 1,\eqno(1.3)
$$
for some constant $c>0$ and some pgf $Q(z)$ satisfying $Q(0)=0$. Also, a distribution with support in {\br} is {\gid} if and only if its Laplace-Stieltjes transform (LST) has the form
$$
\phi(u)=\bigl\{1+\psi(u)\bigr\}^{-1},\qquad u\ge 0,\eqno(1.4)
$$
where $\psi(u)$ has a completely monotone derivative with $\psi(0)=0$.

Let $N_\la(\cdot)$ be a Poisson process of intensity $\la$ and $T$ be an \br-valued rv independent of $N_\la(\cdot)$. The \bz-valued rv $N_\la(T)$ is called a $\la$-Poisson mixture with mixing rv $T$. Its pgf is given by
$$
P_{N_\la(T)}(z)=\phi_T(\la(1-z)),\eqno(1.5)
$$
where $\phi$ is the LST of $T$.

\vskip .5 cm

\goodbreak\noi {\bf 2. Quasi-geometric infinite divisibility}

\smallskip

\Def 2.1. Let $\la>0$. An \br-valued rv $X$ is said to have a $\la$-quasi-geometric infinitely divisible (\lqgid) distribution if the corresponding $\la$-Poisson mixture $N_\la(X)$ is {\gid}

Using (1.4) it can be easily shown that a distribution on {\br} with LST $\phi$ is {\gid} if and only if for any $\tau>0$,
$$
(-1)^nK^{(n)}(\tau)\ge 0, \qquad n\ge 0 \eqno (2.1)
$$
where $K(\tau)=K^{(0)}(\tau)=-{\phi'(\tau)\over\phi(\tau)^2}$ and $K^{(n)}(\tau)$ is its $n$-th derivative, $n\ge 1$. It turns out that (2.1) with $\tau=\la$ characterizes the {\lqgid} property as the following result shows.

\Prop 2.2. Let $X$ be an \br-valued rv with LST $\phi(\tau)$ and $\la>0$. The following assertions are equivalent.

\noi (i) $X$ has a {\lqgid} distribution;

\noi (ii) Condition (2.1) holds for $\tau=\la$ and hence for any $0<\tau\le \la$.

\noi (iii) For any $0<p<1$, the function 
$$
G(z;\la,p)={\phi(\la(1-z))\over p+q\phi(\la(1-z))},\qquad |z|\le 1, \eqno(2.2)
$$ 
is a pgf, where $q=1-p$.

\noi (iv) $N_\la(X)$ satisfies the stability equation 
$$
N_\la(X)\buildrel d \over =B_\la\,(N_\la(X)+S_\la),
$$
for some \bz-valued rv $S_\la$ and some mixed Bernoulli variable $B_\la$ with mixing variable $W_\la$ taking values in $(0,1)$ and with mean $0<E(W_\la)<1$. The rv's $N_\la(X)$, $B_\la$, and $S_\la$ are assumed independent.

\proof (i)$\Leftrightarrow$(ii) By (1.3) and (1.5), the Poisson mixture $N_\la(X)$ is {\gid} if and only if its pgf satisfies 
$$
\phi(\la(1-z))=\bigl(1+c_\la(1-Q_\la(z)\bigr)^{-1},\qquad 0\le z\le1,\eqno(2.3)
$$
for some constant $c_\la>0$ and some pgf $Q_\la(z)$ such that $Q_\la(0)=0$. Assume that $X$ is {\lqgid}. Solving for $Q_\la(z)$ in (2.3) yields
$$
Q_\la(z)=1+c_\la^{-1}-\bigl(c_\la\phi(\la(1-z))\bigr)^{-1}.\eqno(2.4)
$$
Direct calculations imply that the higher order derivatives of $Q_\la(z)$ are given by
$$
Q_\la^{(n+1)}(z)=(-1)^{n}\la^{n+1}c_\la^{-1} K^{(n)}(\la(1-z)),\qquad \ n\ge 0,\eqno(2.5)
$$
and hence, since $Q_\la(z)$ is a pgf,
$$
(-1)^{n}K^{(n)}(\la(1-z))\ge 0,
$$
for any $0\le z<1$ and any $n\ge 0$. Now any $0<\tau\le \la$ can be written as $\tau=\la(1-z)$ for some $z$ ($z=1-{\tau\over \la}$). Therefore (ii) follows. Conversely, assume (2.1) holds for $\tau=\la$. In view of the fact that $0\le 1-\phi(\la(1-z))< 1$ for $0\le z\le 1$ and that $\phi(\la(1-z))$ is a pgf, $Q_\la(z)$ of (2.4), with $c_\la={1\over \phi(\la)}-1$, admits a power series expansion whose coefficients are necessarily nonnegative by (2.5) and (2.1) applied at $z=0$ and $\tau=\la$ repectively. This implies that $Q_\la(z)$ is itself a pgf. 

\noi (i)$\Leftrightarrow$(iii) It is easy to see that (2.2) is equivalent to 
$$
\phi(\la(1-z))={pG(z;\la,p) \over 1-qG(z;\la,p)}, \qquad |z|\le 1.\eqno(2.6)
$$
By definition, $X$ is {\lqgid} if and only if for any $0<p<1$, $G(z;\la,p)$ is the pgf of the rv $X_i^{(p)}$ in (1.1). Finally, (i)$\Leftrightarrow$(iv) follows from Proposition 2.1. in Aly and Bouzar (2000). \endproof

Following Goldie and Puri (1979), we denote by $\calgl$ the class of {\lqgid} distributions on {\br}. We also let
$$
\calgstar=\bigcup_{\la>0}\calgl,\quad \hbox{and}\quad \calginf=\bigcap_{\la>0}\calgl.\eqno(2.7)
$$

\Cor 2.3. (i) For any $0<\la_1<\la_2$, $\calglt\subset\calglo$.

\noi (ii) $\calginf$ is the set of all {\gid} distributions on {\br}.

\noi(iii) A distribution on {\br} is {\gid} if and only if it is {\lqgid} for every $\la$ in an unbounded subset of $(0,\infty)$.

\noi (iv) Let $\la>0$. A distribution on {\br} with LST $\phi(\cdot)$ is {\gid} if and only if for every $p\in (0,1)$, $G(z;\la,p)$ of (2.2) is the pgf of a $\la$-Poisson mixture. In this case the mixing distribution is itself {\gid} with LST $\phi_p(\tau)={\phi(\tau)\over p+q\phi(\tau)}$.

\proof Part (i) follows from Proposition 2.2((i)$\Leftrightarrow$(ii)). Since a distribution on {\br} is {\gid} if and only if (2.1) holds for any $\tau >0$, Proposition 2.2 implies (ii). We note that by (i), $\displaystyle\calginf=\bigcap_{\la\in A}\calgl$ for any unbnounded subset of $A$ of $(0,\infty)$, and thus (iii) is equivalent to (ii). Finally, to prove (iv), if $\phi(\tau)$ is {\gid}, then by (2.2) and (1.4) $\phi_p(\tau)={\phi(\tau)\over p+q\phi(\tau)}={1\over 1+p\psi(\tau)}$ and $G(z;\la,p)={1\over 1+p\psi(\la(1-z))}$, where $\psi(\tau)$ has a completely monotone derivative (and $\psi(0)=0$). Therefore, again by (1.4), $\phi_p(\tau)$ is the LST of a {\gid} distribution and hence $G(z;\la,p)$ is the pgf of a $\la$-Poisson mixture. Conversely, assume that for any $0<p<1$, $G(z;\la,p)$ is the pgf of a $\la$-Poisson mixture, then by Lemma A.6 in van Harn and Steutel (1993), $\phi_p(\tau)={\phi(\tau)\over p+q\phi(\tau)}$, $\tau\ge 0$, $q=1-p$, is the LST of a distribution on {\br}. Therefore, by (1.1), $\phi(\tau)$ is the LST of a {\gid} distribution (with $\phi_p(\tau)$ being the LST of $X_i^{(p)}$) .\endproof

Contrasting Proposition 2.2 (iii) and Corollary 2.3 (iv), it is worth remarking that if for some $\la>0$ a distribution on {\br} with LST $\phi(\cdot)$ is {\lqgid} but not {\gid}, then there must exist $0<p<1$ such that the pgf $G(z;\la,p)$ of (2.2) is not a $\la$-Poisson mixture.

Puri and Goldie (1979) obtained necessary an sufficient conditions for an i.d. discrete distribution on {\bz} to be a Poisson mixture generated by an i.d. mixing distribution. We state an analogous result for discrete {\gid} distributions.

\Prop 2.4. Let $P(z)$ be a pgf. Then $P(z)$ is the pgf of a Poisson mixture generated by a {\gid} mixing distribution with LST $\phi(\tau)$ if and only if the two conditions below hold:

\noi (i) $P(z)$ is defined and satisfies $0<P(z)\le 1$ for all $z\in (-\infty,1]$;

\noi (ii) the mapping $H(z)={1\over P(z)}-1$ is in $C^\infty((-\infty,1))$ and
$$
H^{(n)}(z)\le 0, \quad \hbox{for all } n\ge 1,\ \hbox{and}\  \hbox{all } z\in(-\infty,1).\eqno(2.8)
$$
In this case $P(z)$ is necessarily {\gid} Moreover, for any $p\in (0,1)$, the pgf $G_p(z)={P(z)\over p+qP(z)}$ (where $q=1-p$) is also a Poisson mixture generated by a {\gid} mixing distribution with LST $\phi_p(\tau)={\phi(\tau)\over p+q\phi(\tau)}$.

\proof Suppose $P(z)$ is the pgf of a Poisson mixture generated by a {\gid} mixing distribution with LST $\phi(\tau)$. Then $P(z)=\phi(\la(1-z))$ for some $\la >0$ from which (i) and the first part of (ii) follow trivially. By (1.4), $\phi(\tau)={1\over 1+\psi(\tau)}$ where $\psi(\tau)$ has a completely monotone derivative on $[0,\infty)$ (and $\psi(0)=0$). Hence $H(z)=1/P(z)=1+\psi(\la(1-z))$, $z\in (-\infty,1)$ which implies that for any $n\ge 1$ and $z\in (-\infty, 1)$, 
$$
H^{(n)}(z)=\la^n(-1)^n\psi^{(n)}(\la(1-z)),\eqno(2.9)
$$
which in turn implies (2.8). The fact that $P(z)$ is itself {\gid} is a consequence of Proposition 4.2 in Aly and Bouzar (2000). It easily follows that $G_p(z)$ is a Poisson mixture whose mixing distribution has LST $\phi_p(\tau)$ and is necessarily {\gid} (as seen in the proof of Corollary 2.3(iv)). Conversely, suppose that (i) and (ii) hold. Define $\phi(\tau)=P(1-\tau)$ and $\psi(\tau)={1\over \phi(\tau)}-1=H(1-\tau)$ for $\tau\ge 0$. Then $P(z)={1\over 1+\psi(1-z)}$. By Proposition 4.4 in Aly and Bouzar (2000), it is sufficient to prove that $\psi(\tau)$ has a completely monotone derivative. Trivially, for any $n\ge 0$ and $\tau>0$,
$$
(-1)^n\psi^{(n+1)}(\tau)=-H^{(n+1)}(1-\tau)\ge 0,
$$
where the latter inequality follows from (2.8).\endproof

Next, we present some closure properties of $\calgl$ and $\calgstar$ similar to the ones obtained by Puri and Goldie (1979) in the q.i.d. case. In what follows we will say that an \br-valued rv is in ${\cal G}_{(\cdot)}$ if and only if its distribution is in ${\cal G}_{(\cdot)}$.

\Prop 2.5. (i) Let $X$ be an \br-valued rv and $\la>0$. If $X \in \calgl$, then for any positive constant $c$, $cX\in {\cal G}_{\la/c}$. Therefore, for any $0<c<1$, $cX \in \calgl$. Also, $X \in \calgl$ if and only if $\la X\in {\cal G}_1$.

\noi (ii) $\calgstar$ is closed under multiplication by a positive scalar.

\noi (iii) For every $\la>0$, $\calgl$ is closed under convergence in distribution. However, $\calgstar$ is not closed under the same operation.

\proof (i) Let $\phi(\tau)$ be the LST of $X$. The LST of $cX$ for $c>0$ is $\phi_c(\tau)=\phi(c\tau)$. Letting $K_c(\tau)=-{\phi'_c(\tau)\over\phi_c^2(\tau)}$, it follows by (2.1) and Proposition 2.2 that if $X$ is in $\calgl$ for some $\la>0$, then for any $n\ge 0$, $(-1)^n K_c^{(n)}(\la/c)=(-1)^nc^{n+1}K^{(n)}(\la)\ge 0$ which implies that $cX\in {\cal G}_{\la/c}$. The remaining assertions follow from the first part and Proposition 2.2. (ii) is a straightforward consequence of (i). To prove (iii), let $\la>0$ and let $(X_n,n\ge 0)$ be a sequence of \br-valued rv's in $\calgl$ such that $X_n \buildrel d \over \rightarrow X$ for some \br-valued rv $X$. Then by Theorem 10 in Puri and Goldie (1979), $N_\la(X_n) \buildrel d \over \rightarrow N_\la(X)$. Since for every $n\ge 0$, $N_\la(X_n)$ is {\gid}, then by (1.3) and Theorem 1 in Klebanov et al. (1984) (adapted to pgf's), $N_\la(X)$ must be {\gid}. The counterexample for the second part of (iii) is given below as Example (6). \endproof

We conclude this section by giving several examples and counterexamples.

\noi 1) A pgf P(z) which is {\gid} but is not a Poisson mixture: Let $Q(z)$ be the pgf of distribution with support on the nonnegative even integers. Then $P(z)=(1+c(1-Q(z)))^{-1}$, for some $c>0$, is the pgf of a {\gid} distribution but, by Proposition 2.4. (i), it is not a Poisson mixture.

\noi 2) A Poisson mixture with pgf $P(z)$ generated by a {\gid} mixing distribution: $P(z)=(1+c(1-z)^\a)^{-1}$, $0<\a\le 1$ is the pgf of the discrete Mittag-Leffler distribution (Pillai and Jayakumar (1995)) and it is a Poisson mixture generated by a {\gid}, continuous Mittag-Leffler mixing distribution with LST $\phi(\tau)=(1+a\tau^\a)^{-1}$, for some $a>0$ (Pillai (1990), Aly and Bouzar (2000)).

\noi 3) A Poisson mixture generated by a mixing distribution that is in $\calgstar$ but not in $\calginf$: Consider the distribution function $F(x)$ on {\br} with LST  $\phi_\a(\tau)=(1+ \tau^\a)^{-2}$, for some $1/2<\a<1$. Again, this distribution is of the continuous Mittag-Leffler type and hence belongs to $\calginf$. Letting $K_\a(\tau)=-{\phi'_\a(\tau)\over\phi_a^2(\tau)}$, we have $K_\a(\tau)=2\a\tau^{\a-1}(1+\tau^\a)$ and
$$
(-1)^nK_\a^{(n)}(\tau)=2\a\tau^{\a-n-1}(A_n-B_n\tau^\a), \qquad n\ge 1,
$$
where $A_n=\prod_{i=1}^n(i-\a)$ and $B_n=(2\a-1)\prod_{i=2}^n(i-2\a)$. Note that since $1/2<\a<1$, $A_n>0$, and $B_n>0$ for any $n\ge 1$. If $\la$ satifies $0<\la^\a\le{1-\a\over 2\a-1}$, (2.1) holds for $K_\a(\tau)$ at $\tau=\la$, but fails to hold (at $n=1$) for any $\la$ such that $\la^\a>{1-\a\over 2\a-1}$. Hence $F(x)$ belongs to $\calgstar$ but not to $\calginf$.
 
\noi 4) Neither $\calgl$ nor $\calgstar$ is closed under translation: Let $X$ be an \br-valued rv with LST $\phi(\tau)=(1+\log(1+\tau))^{-1}$. It is easy to see that the distribution of $X$ is in $\calginf$. The LST of $X+1$ is $\phi_1(\tau)=e^{-\tau}\phi(\tau)$. Letting $K_1(\tau)=-{\phi'_1(\tau)\over\phi_1^2(\tau)}$, we have $K_1(\tau)=e^\tau\bigr(1+\log(1+\tau)+(1+\tau)^{-1}\bigl)$ and hence, $(-1)K_1'(\tau)=-e^\tau\bigr(1+\log(1+\tau)+(2\tau+1)(1+\tau)^{-2}\bigl)<0$ for all $\tau\ge0$. By Proposition 2.2, $X$ cannot belong to $\calgstar$.

\noi 5) Neither $\calgl$ nor $\calgstar$ is closed under convolution: Let $X$ and $Y$ be \br-valued iid exponentially distributed rv's with mean 1 and common LST $\phi(\tau)=(1+ \tau)^{-1}$. Obviously, the exponential distribution is in $\calginf$. The LST of $X+Y$ is $\phi_1(\tau)=(1+ \tau)^{-2}$ and $K_1(\tau)=-{\phi'_1(\tau)\over\phi_1^2(\tau)}=2(1+\tau)$. Since $(-1)K'_1(\tau)=-2$ for any $\tau>0$, by Proposition 2.2 $X+Y$ cannot belong to $\calgstar$.  

\noi 6) $\calgstar$ is not closed under convergence in distribution: Let $(X_n, n\ge 1)$ be a sequence of \br-valued rv's such that for each $n\ge 1$, $X_n$ has LST $\phi_{\a_n}(\tau)=(1+ \tau^{\a_n})^{-2}$, with $1/2<\a_n<1$ and $\lim_{n\to \infty}\a_n=1$. By Example 3 above, $X_n$ is in $\calgstar$ for every $n\ge 1$. Since $\lim_{n\to \infty}\phi_{\a_n}(\tau)=(1+\tau)^{-2}$, $X_n$ converges in distribution, but its limit is not in $\calgstar$ (see Example 5 above).
\vskip .5 cm

\goodbreak\noi {\bf 3. Log-convexity, log-concavity and quasi-geometric infinite divisibility}

We need the following lemma for our next characterization of the {\lqgid} property. We recall that a sequence of nonnegative real numbers $(b_n, n\ge 0)$ is said to be log-convex (resp. log-concave) if for every $n\ge 1$,
$$
b_{n-1}b_{n+1}\ge b_n^2, \quad (\hbox{resp. } b_{n-1}b_{n+1}\le b_n^2.) \eqno(3.1)
$$
\Lemma 3.1. A distribution $(p_n, n\ge 0)$ on {\bz}, $0<p_0<1$, is {\gid} if and only if the sequence $(a_n, n\ge 0)$ defined by
$$
p_{n+1}=\sum_{k=0}^np_ka_{n-k}, \quad n\ge 0, \eqno(3.2)
$$
is nonnegative and necessarily satisfies $\displaystyle \sum_{n=0}^\infty a_n<\infty$. Consequently, if $(p_n, n\ge 0)$ is log-convex, then $(p_n, n\ge 0)$ is {\gid}

\proof By (1.3), $(p_n, n\ge 0)$ is {\gid} if and only if its pgf $P(z)$ satisfies 
$$
(1+c)P(z)-cP(z)Q(z)=1, \qquad |z|\le 1,\eqno(3.3)
$$
for some pgf $Q(z)$ with distribution $(q_n, n\ge 0)$, $q_0=0$ and some constant $c>0$. Using the power series expansions of $P(z)$ and $Q(z)$ in (3.3) yields (3.2) with $a_n={c\over 1+c}q_{n+1}$, $n\ge 0$. By passing to generating functions in (3.2) we obtain (3.3) and hence the converse. The second part follows from Lemma 4.2.2 in Steutel (1970).\endproof

We denote by $(p_n(\tau), n\ge 0)$ the distribution of a $\tau$-Poisson mixture, $\la>0$, generated by the mixing distribution $F(x)$ on {\br}, i.e.,
$$
p_n(\tau)={1\over n!}\int_0^\infty (\tau x)^n e^{-\tau x}\,dF(x), \quad n\ge 0. \eqno(3.4)
$$

\Prop 3.2. A distribution function $F$ with support in {\br} is {\lqgid} if and only if the sequence $(a_n(\tau), n\ge 0)$, $\tau>0$, defined by (3.2) (with $p_n=p_n(\tau)$) is nonnegative for $\tau=\la$ and hence for every $0<\tau\le \la$. In particular, if $(p_n(\la), n\ge 0)$ is log-convex, then $F$ is {\lqgid}.

\proof Straightforward from Proposition 2.2 and Lemma 3.1. \endproof

We recall that a distribution $(p_n, n\ge 0)$ on {\bz} is infinitely divisible (i.d.) if and only if the sequence $(r_n, n\ge 0)$ defined by
$$
(n+1)p_{n+1}=\sum_{k=0}^nr_kp_{n-k}, \quad n\ge 0, \eqno(3.5)
$$
is nonnegative (see for example Steutel (1970)) and satisfies necessarily $\sum_{n=0}^\infty r_n(n+1)^{-1}<\infty$. We will refer to the sequence $(r_n, n\ge 0)$ as the L\'evy measure of (an i.d. distribution) $(p_n, n\ge 0)$.  

\Lemma 3.3. Let $(p_n, n\ge 0)$ be an i.d. distribution on {\bz} and $(r_n, n\ge 0)$ its associated L\'evy measure. Then $(p_n, n\ge 0)$ is {\gid} if and only if the sequence $(b_n,n\ge 0)$ defined by $b_0=0$ and  
$$
(n+1)b_{n+1}=r_n-\sum_{k=1}^nb_kr_{n-k}, \quad n\ge 0,\eqno (3.6)
$$
is nonnegative and necessarily satisfies $\displaystyle \sum_{n=0}^\infty b_n<1$.

\proof Again, by (1.3) $(p_n,n\ge 0)$ is {\gid} if and only if its pgf $P(z)$ satisfies
$$
(1+c(1-Q(z)){d \over d\,z} \ln P(z)=cQ'(z),\quad |z|<1,\eqno(3.7)
$$
where $Q(z)$ is the pgf of a distribution $(q_n,n\ge 0)$ on {\bz}, $Q(0)=0$, and $c>0$. Noting that $(r_n, n\ge 0)$ is the sequence of the coefficients of the power series expansion of ${d \over d\,z} \ln P(z)$, $|z|<1$, it can be easily deduced that (3.6) and (3.7) are equivalent (by letting $b_n={c\over 1+c}q_n$) via power series representations. Clearly, $\sum_{n=0}^\infty b_n={c\over 1+c}<1$. \endproof

\Prop 3.4. Let $F$ be a $\la$-q.i.d. distribution function on {\br}. Let $(p_n(\la), n\ge 0)$ be the corresponding Poisson mixture, as given by (3.4), and $(r_n(\la),n\ge 0)$ its L\'evy measure. $F$ is {\lqgid} if and only if the sequence $(b_n(\tau), n\ge 0)$, $\tau>0$, defined by (3.6) (with $r_n=r_n(\tau)$) is nonnegative for $\tau=\la$, and hence for any $0<\tau\le \la$.

\proof Straightforwardly form Proposition 2.2 and Lemma 3.3.\endproof
 
\Lemma 3.5 Let $(p_n, n\ge 0)$ be an i.d. distribution on {\bz} and $(r_n, n\ge 0)$ its associated L\'evy measure. Assume that $(r_n, n\ge 0)$ is log-convex. The following assertions are equivalent.

\noi (i)  $(p_n, n\ge 0)$ is {\gid};

\noi (ii) $r_0^2\le r_1$;

\noi (iii) $(p_n, n\ge 0)$ is log-convex.

\proof (ii)$\Leftrightarrow$(iii) is simply Theorem 2 in Hansen (1988). We prove (i)$\Leftrightarrow$(ii). If $(p_n, n\ge 0)$ is {\gid}, then by (3.6) applied to $n=0,1$, $b_1=r_0$ and $2b_2=r_1-b_1r_0$. Since $b_2\ge 0$, it follows that $r_0^2\le r_1$. Conversely, assume that $(r_n, n\ge 0)$ is log-convex and that $r_0^2\le r_1$. By Lemma 3.3, it is sufficient to prove that $(b_n, n\ge 0)$ of (3.6), with $b_0=0$, is nonnegative. We proceed by induction. We have trivially, $b_0, b_1\ge 0$. Assume $b_k\ge 0$, $0\le k\le n$. For $n\ge 0$, let $A_n=\prod_{i=0}^nr_i$. Then by (3.6)
$$
(n+1)b_{n+1}A_{n-2}=A_{n-2}r_n-\sum_{k=1}^nb_kA_{n-2}r_{n-k}, \quad n\ge 2.\eqno(3.8)
$$
By applying (3.1) to $r_n$, $n\ge 2$, and letting $A_{-1}=0$,
$$
A_{n-2}r_n=A_{n-3}r_{n-2}r_n\ge A_{n-3}r_{n-1}^2.\eqno(3.9)
$$
Applying repeatedly (3.1) to $r_{(\cdot)}$ yields for any $1\le k\le n-1$,
$$
\eqalign{
A_{n-2}r_{n-k}&=(r_0\cdots r_{n-k-1}r_{n-k+1}\cdots r_{n-2})r_{n-k}^2\cr
&\le (r_0\cdots r_{n-k-1})r_{n-k+1}\cdots r_{n-2}r_{n-k-1}r_{n-k+1}\cr
&=(r_0\cdots r_{n-k-1})r_{n-k+1}^2 r_{n-k+2}\cdots  r_{n-2}r_{n-k-1}\cr
&\le (r_0\cdots r_{n-k})r_{n-k+2}^2 r_{n-k+3}\cdots  r_{n-2}r_{n-k-1}\cr
&\le \cdots\cdots\cdots\le A_{n-3}r_{n-1}r_{n-k-1}.
}\eqno(3.10)
$$
The case $k=n$ follows similarly via the inequalities $r_0^2\le r_1$ and (3.1):
$$
\eqalign{
A_{n-2}r_0&=(r_1\cdots r_{n-2})r_0^2\le (r_2\cdots r_{n-2})r_1^2\cr
&\le (r_2r_3\cdots r_{n-2})r_0r_2=(r_0r_3\cdots r_{n-2})r_2^2\cr
&\le(r_0r_3r_4\cdots r_{n-2})r_1r_3\le \cdots\cdots\le A_{n-3}r_{n-1}\cr
}
\eqno(3.11)
$$
Hence by (3.8)--(3.11) and the induction hypothesis,
$$
\eqalign{
(n+1)b_{n+1}A_{n-2}& \ge A_{n-3}r_{n-1}^2-\sum_{k=1}^{n-1}b_kA_{n-3}r_{n-1}r_{n-k-1}-b_nA_{n-3}r_{n-1}\cr
& \ge A_{n-3}r_{n-1}[(r_{n-1}-\sum_{k=1}^{n-1}b_k r_{n-k-1})-b_n]\cr
&=A_{n-3}r_{n-1}(n-1)b_n,\cr
}
$$
where the last equation follows from (3.6). Since $A_k$'s and the $r_k$'s are nonnegative, it follows that $b_{n+1}\ge 0$. \endproof

\Prop 3.6. Let $F$ be a $\la$-q.i.d. distribution function on {\br}. Let $(p_n(\la), n\ge 0)$ be the corresponding Poisson mixture, as given by (3.4), and $(r_n(\la),n\ge 0)$ its L\'evy measure. Moreover, assume that $(r_n(\la),n\ge 0)$ is log-convex. The following assertions are equivalent.

\noi(i) $F$ is {\lqgid};

\noi(ii) $r_0^2(\la)\le r_1(\la)$;

\noi(iii) $(p_n(\la), n\ge 0)$ is log-convex.

\proof Again, this follows straightforwardly form Proposition 2.2 and Lemma 3.5.\endproof

The log-concave versions of Lemma 3.5 and Lemma 3.7 are established similarly to their log-convex counterparts. We state the results without proofs.

\Lemma 3.7. Let $(p_n, n\ge 0)$ be an i.d. distribution on {\bz} and $(r_n, n\ge 0)$ its associated L\'evy measure. Assume that $(r_n, n\ge 0)$ is log-concave. The following assertions are equivalent.

\noi (i)  $(p_n, n\ge 0)$ is {\gid};

\noi (ii) $r_0^2\ge r_1$;

\noi (iii) $(p_n, n\ge 0)$ is log-concave.
 
\Prop 3.8. Let $F$ be a $\la$-q.i.d. distribution function on {\br}. Let $(p_n(\la), n\ge 0)$ be the corresponding Poisson mixture, as given by (3.4), and $(r_n(\la),n\ge 0)$ its L\'evy measure. Moreover, assume that $(r_n(\la),n\ge 0)$ is log-concave. The following assertions are equivalent.

\noi(i) $F$ is {\lqgid};

\noi(ii) $r_0^2(\la)\ge r_1(\la)$;

\noi(iii) $(p_n(\la), n\ge 0)$ is log-concave.

\Remark Lemma 3.5 and Lemma 3.7 strengthen Theorem 1 and Theorem 2 obtained by Hansen (1988).

\vskip .5 cm

\noi {\bf 4. Generalizations}\par\nobreak

\smallskip

A more general notion of infinite divisibility based on (1.1) was studied by several authors (see Gnedenko and Korolev (1996), Section 4.6, for details and further references.) The definition is as follows. Let $I$ be a subset of $(0,1)$ and let ${\cal N}=\{N_p, p \in I\}$ be a family of \bz-valued rv's such that $E(N_p)=1/p$ for any $p\in I$,  and 
$$
H_{p_1}\circ H_{p_2}(z)=H_{p_2}\circ H_{p_1}(z),\quad \hbox{for any $p_1,\ p_2 \in I$}, \eqno(4.1)   
$$
where $H_p$ is the pgf of $N_p$. A rv $X$ is said to be $\cal N$-infinitely divisible, or {\nid}, if it satisfies (1.1) for any $N_p\in \cal N$. We recall that (4.1) implies (see the proof of Theorem 4.6.1 in Gnedenko and Korolev (1996)) the existence of an LST $\var$ satisfying $\var(0)=-\var'(0)=1$ and 
$$
\var(\tau)=H_p(\var(p\tau)), \quad \hbox{for any $\tau >0$ and $p\in I$.} \eqno(4.2)
$$
Aly and Bouzar (2000) showed that a {\bz} (resp. {\br})-valued rv $X$ with pgf $P(z)$ (resp. LST $\phi(\tau)$) is {\nid} if and only if
$$
P(z)=\var[c(1-Q(z))] \quad (\hbox{resp. } \phi(\tau)= \var[\psi(\tau)]), \eqno(4.3)
$$
where $\var(\cdot)$ is as in (4.2) and $Q(z)$ (resp. $\psi(\tau)$) is a pgf satisfying $Q(0)=0$ and $c>0$ (resp. has a completely monotone derivative with $\psi(0)=0$).

Similarly to the geometric case, an \br-valued rv $X$ is said to have a {\lqnid} distribution for $\la>0$ if the corresponding $\la$-Poisson mixture $N_\la(X)$ is {\nid} 

Next, we state several characterizations of the {\nid} property. The proof follows from (4.3) and the arguments used in the proof of Proposition 2.2. The details are omitted. 

\Prop 4.1. Let $X$ be an \br-valued $X$ with LST $\phi(\tau)$ and $\la>0$. The following assertions are equivalent.

\noi (i) $X$ has a {\lqnid} distribution;

\noi (ii) Condition (2.1) holds for $\tau=\la$ and hence for any $0<\tau\le \la$, where
$$
K(\tau)={\phi'(\tau)\over \var'[\var^{-1}(\phi(\tau))]};\eqno(4.4)
$$

\noi (iii) For any $p\in I$, the function
$$
G(z;\la,p)=H_p^{-1}[\phi(\la(1-z))],\qquad |z|\le 1, \eqno(4.5)
$$ 
is a pgf.

Defining $\calgl$, $\calgstar$, and $\calginf$ for the {\nid} case similarly to their geometric counterparts (see (2.7)), we have the following corollary that can be proven along the same lines as Corollary 2.3, via (4.2). Again, the details are omitted.

\Cor 4.2. (i) For any $0<\la_1<\la_2$, $\calglt\subset\calglo$.

\noi (ii) $\calginf$ is the set of all {\nid} distributions on {\br}.

\noi(iii) A distribution on {\br} is {\nid} if and only if it is {\lqnid} for every $\la$ in an unbounded subset of $(0,\infty)$.

\noi (iv) Let $\la>0$. A distribution on {\br} with LST $\phi(\cdot)$ is {\nid} if and only if for every $p\in I$, $G(z;\la,p)$ of (4.5) is the pgf of a $\la$-Poisson mixture. In this case the mixing distribution is itself {\nid} with LST $\phi_p(\tau)=H_p^{-1}(\phi(\tau))$.

Next, we state without proof the analogue of Proposition 2.4. 

\Prop 4.3. Let $P(z)$ be a pgf. Then $P(z)$ is the pgf of a Poisson mixture generated by an {\nid} mixing distribution with LST $\phi(\tau)$ if and only if the two conditions below hold:

\noi (i) $P(z)$ is defined and satisfies $0<P(z)\le 1$ for all $z\in (-\infty,1]$;

\noi (ii) the mapping $H(z)=\var^{-1}(P(z))$, with $\var(\cdot)$ of (4.2), is in $C^\infty((-\infty,1))$ and
$$
H^{(n)}(z)\le 0, \quad \hbox{for all } n\ge 1,\ \hbox{and}\  \hbox{all } z\in(-\infty,1).\eqno(4.6)
$$
In this case $P(z)$ is necessarily {\nid} Moreover, for any $p\in I$, the pgf $G_p(z)=H_p^{-1}(P(z))$ is also a Poisson mixture generated by an {\nid} mixing distribution with LST $\phi_p(\tau)=H_p^{-1}(\phi(\tau))$.

We also note that the closure properties obtained in Proposition 2.5 generalize verbatim to 
quasi-$\cal N$-infinite divisibility.

We conclude by noting that classical (resp. geometric) infinite divisibility corresponds to the family of rv's $\cal N$ where $N_p={1\over p}$ with probability 1, $p\in I=\{{1\over n}: n\ge 1\}$ (resp. $N_p$ has distribution (1.2), $p\in I=(0,1)$). In the classical case $\var(\tau)=e^{-\tau}$ and the results of these sections thus include those of Puri and Goldie (1979), whereas in the geometric case, $\var(\tau)=(1+\tau)^{-1}$ thus yielding the results of Section 2.

\goodbreak\bibs

\ref E.E. Aly and N. Bouzar, On geometric infinite divisibility and stability, Annals of the Institute of Statistical Mathematics, 52 (2000), 790--799.

\ref B.V. Gnedenko and V.Y. Korolev, {\it Random Summation: Limit Theorems and Applications}, CRC Press, Boca Raton, New York, London, Tokyo (1996).

\ref B.G. Hansen, On log-concave and log-convex infinitely divisible sequences and densities, Annals of Probability, 16 (1988), 1832--1839.

\ref P.A. Jacobs,(1986). First passage times for combinations of random loads, SIAM Journal of Applied Mathematics, 46 (1986), 643--656.

\ref N.L. Johnson, S. Kotz, and A.W. Kemp, {\it Univariate Discrete Distributions}, Second ed., John Wiley and Sons, New York (1992).

\ref K. van Harn and F.W. Steutel, Stability equations for processes with stationary independent increments using branching processes and Poisson mixtures, Stochastic Processes and their Applications, 45 (1993), 209--230.

\ref Y. Kebir, Laplace transforms and the renewal equation, Journal of Applied Probability, 34 (1997), 395--403.

\ref L.B. Klebanov, G.M. Maniya, and I.A. Melamed, A problem of Zolotarev and analogs of infinitely divisible and stable distributions in a scheme for summing a random number of random variables, Theory of Probability and Applications, 29 (1984), 791--794.

\ref T.J. Kozubowski and S.T. Rachev, Univariate geometric stable laws, Journal of Computational Analysis and Applications, 1 (1999), 177--217.

\ref T.J. Kozubowski and S.T. Rachev, The theory of geometric stable distributions and its use in modeling financial data, European Journal of Operations Research, 74 (1994), 310--324.

\ref A.G. Pakes, (1995), Characterization of discrete laws via mixed sums and Markov branching processes, Stochastic Processes and their Applications, 55 (1995), 285--300.

\ref R.N. Pillai, On Mittag-Leffler functions and related distributions, Annals of the Instistute of Statistical Mathematics, 42 (1990), 157--161.

\ref R.N. Pillai and K. Jayakumar, Discrete Mittag-Leffler distributions, Statistics and Probability Letters, 23 (1995), 271--274.

\ref P.S. Puri and C.M. Goldie, Poisson Mixtures and quasi-infinite divisibility of distributions, Journal of Applied Probability, 16 (1979), 138--153.

\ref F.W. Steutel, {\it Preservation of Infinite Divisibility under Mixing}, Mathematical Centre Tract {\bf 33}, Mathematisch Centrum Amsterdam (1970).

\bye